# Research on Failure Characteristics of Circular Expansion Foundation of Large Wind Turbine

## Peng Cheng[1], Zijian Fan[2],*

*1 School of Civil and Engineering, University of South China, Hengyang 421001, China*
*2 School of Civil and Engineering, University of South China, Hengyang 421001, China*

## 1. Summary

Wind power as clean energy has been widely promoted. However, during the operation of a wind turbine in a mountainous area, a large number of instability phenomena appear in the wind turbine foundation. To study the damage characteristics of a mountain wind turbine foundation, the influence of topography on a mountain wind turbine is considered. According to the complex mechanical structure characteristics of the wind turbine, a circular expansion foundation model have a ratio of 1:10 to the wind turbine which 2MW installed capacity is established. The foundation model is simulated by the ABAQUS, and the impact of soil pressure on the wind turbine foundation under random wind load is analyzed. Understand the fragile parts of the wind turbine foundation. The results show that :(1) under the action of random wind load, the strain of the anchor rod, and the ground ring, of the bottom plate is small and will not be destroyed ;(2) under the influence of wind load, he soil within a 1.5 times radius range of wind turbine foundation is easily disturbed, but the core soil below the foundation is unaffected. It is suggested that the soil within a 1.5 times radius range of the wind turbine foundation should be strengthened to ensure the safety and stability of the foundation.

## 2. Introduction

The wind turbine is a power generation facility built in a resource-rich area of the wind where the wind loads are highly variable. The direction and strength of the wind force are always fluctuating over time which makes hard to simulate the force of each component of the wind turbine, A change in the angle of the wind from 30 degrees to 40 degrees will result in a 27% change in the thrust experienced by the wind turbine[1-2].And the wind at different altitudes is different[3-5].So that the wind turbine structure as a force-complicated structure is easy to damage during work. Especially the wind power foundation as the main force structure, not only to withstand large vertical loads but also to withstand large horizontal force and tipping moment[6]. Because of this characteristic, the wind turbine foundation and the part of the foundation connected to the tower that are highly vulnerable to damage during working process.

To study the dynamic response characteristics of wind turbine structure, a lot of of research has been done that mainly focused on the different parts of the wind turbine[7-9]. However, domestic and foreign scholars had carried out some researches on the basic destruction characteristics of wind turbines. The failure mode of wind turbine wide and shallow cylinder foundation under bending moment load was studied[10-12]. Xu et al.[13] analyzed the dynamic response of three different forms of offshore wind turbine pile foundation under seismic loading. Wang et al.[14] analyzed the development and failure modes of the surrounding soil of the single pile foundation in the soft soil multi-layer foundation under different loads. Wu et al.[15] focus on the failure model about offshore wind turbine single pile foundation which under soft soil condition and cyclic load. Wu[16] analyzed the damage process and mechanism of the onshore wind turbine coupling system foundation. Lu et al.[17] and He et al.[18] analyzed fatigue damage effect about foundation ring and concrete joint based on project example. Stamatopoulos[19]investigated the seismic behavior of wind turbine supported by a circular spread foundation using the current Greek aseismic code. Dai[20] evaluated the structural responses of the embedded ring exhibited significant vertical displacements during periods of turbulent wind speed.

Chinese research on wind-base damage characteristics was mostly focused on offshore wind foundations and many theoretical studies and numerical simulations based on specific working conditions. Research on the onshore wind turbine foundation was far less than the research on the foundation of the offshore wind turbine foundation. And the broken model is verity in different working conditions[21-22]. Like-wise, the mountain wind turbine foundation type is various, but the most popular one is the expansion foundation[3].

In this paper, the circular expansion foundation of 2MW installed capacity was taken as the research object, and a 1:10 scale model was established. Firstly, focusing on the foundation and connection position between foundation and tower, the bolt, foundation ring and base plate strain, base displacement and base pressure were tested in all aspects. Then, the earth pressure changes were analyzed by the numerical model. Finally, the failure mechanism of the circular expansion foundation of large wind turbines is summarized.

# 3. Materials and Methods

## 3.1 Simplification of load

The wind foundation is affected by two kinds of loads, permanent load and variable load. The permanent load include foundation gravity, superstructure vertical load and backfill gravity. They were represented by G1, Fz and G2. The variable load include loads in the x and y directions caused by wind loads, moments in the x, y and z directions. They were represented by Fx, Fy, Mx, My and Mz. The representation of the six loads other than the two gravitational forces in the basic coordinate system of the wind turbine was shown in Figure 1. According to Chinese wind turbine foundation design requirements, and because the two vertical load could combine as one load and the two vertical moments could combine as one moment when sign the wind turbine foundation, then, the load was simplified to a vertical load, a horizontal force and a horizontal moment. They were represented by Fz, Fr and Mr.

Since the value of Mz was generally an order of magnitude smaller than Mr, according to China's number FD 003-2007 wind turbine foundation design requirements, it is usually not considered. In the foundation design, the two horizontal forces (Fx and Fy) and the horizontal moments (Mx and My) of the same working condition could be combined into the horizontal resultant force Fr and the horizontal combined torque Mr. Thus, the load acting on the foundation of the wind turbine could be reduced to a vertical force Fz, a horizontal resultant force Fr and a horizontal combined torque Mr, in addition to the self-weight.

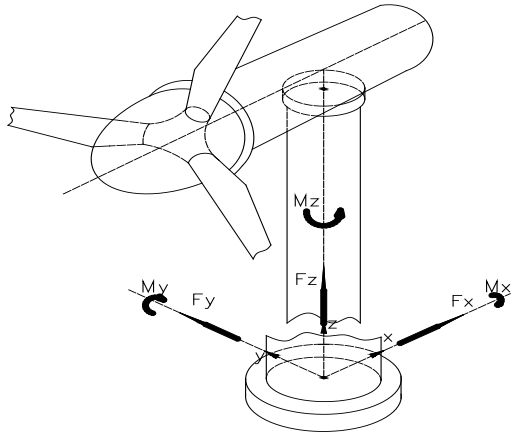

Figure 1. The load acting on the foundation

The calculation of the horizontal thrust of the wind turbine and the wind load of the tower by hydrodynamics calculation software and theoretical methods. It was found that the wind load acting on the tower is usually less than 20% of the horizontal wind load of the entire wind turbine. The amount of change in the wind load of the tower is small relative to the thrust of the wind turbine, so it could be assumed that the wind load of the tower was also static. Therefore, Fr was regarded as the dynamic load that changes according to the horizontal thrust of the wind wheel in the subsequent study of this paper. And Fz, Mr, G1 and G2 were regarded as static loads acting on the foundation of the wind turbine.

To accurately simulate the wind load conditions in the natural state, and highlight the additional effect of the wind load on the foundation. It was assumed that the horizontal wind load was a pulsating wind curve that considers both the mean wind and the random fluctuation of the periodic variation, and the sinusoid with a variation period of 10 min and a variation amplitude of 0.2 Fr. After used MATLAB software simulation, the wind force curve is shown in Figure 2.

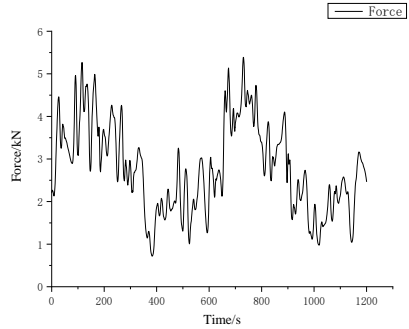

Figure 2. The level random wind load

## 3.2 Test material

In this test, the foundation soil was made of undisturbed slope soil filtered through a 2 mm sieve hole, and a certain amount of sand was added to it to increase the density of the foundation soil. The soil gradation curve used in the test is shown in Figure 3. The soil particle non-uniformity coefficient is 10.2.

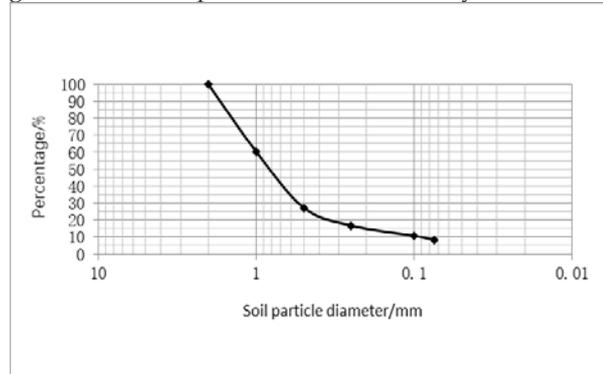

Figure 3. Particle grading curve

The test used a model in which the elastic modulus and the geometric parameters are completely similar and the density is similar. All geometric parameters and material parameters of the model are shown in Table 1.

Table 1. Model material mechanics parameters

| Position | Material | Parameter | Numerical value |
|---|---|---|---|
| Tower | Steel | E | 212.8GPa |
| | | $\mu$ | 0.31 |
| | | $\sigma_s$ | 349.2MPa |
| | | $\varrho$ | 7.85g·cm-3 |
| Reinforced concrete foundation | Concrete | E | 32.3 MPa |
| | | $\mu$ | 0.23 |
| | | $f_{cu.k}$ | 27.7 MPa |
| | | $\varrho$ | 0.0234 g·m-3 |
| | Reinforcement | E | 203.8 MPa |
| | | $\mu$ | 0.31 |
| | | $\sigma_s$ | 237.6 MPa |
| | | $\varrho$ | 7.85 g·cm-3 |
| Foundation soil | Undisturbed soil | $E_s$ | 8.35 MPa |
| | | $\mu$ | 0.33 |

| | |
|---|---|
| ψ | 23.2 |
| ϱ | 2.03 g·cm-3 |
| c | 17.6 kPa |
| ω | 16.7% |

To highlight the main part of this test and simplify the experimental conditions. The tower ignored the ancillary facilities on the tower and only tests the first and second sections of the bottom tower, and included connecting parts such as foundation rings and bolts. And consider the wall thickness is only 3mm after the scale, To ensure the similarity of its parameters, the steel used the same material as the prototype tower. The reinforced concrete foundation was simulated by simulated concrete, and the largest aggregate particle size was selected as 1/10 of the actual maximum aggregate particle size. The plain concrete of the cushion layer was simulated by cement mortar, and the steel wire for steel reinforcement was used to ensure the reinforcement ratio. In addition to optimizing the gradation of the soil, the simulation of the foundation soil also tampered with the soil by artificial tamping to ensure the compactness of the soil.

### 3.3 Loading device and observation instrument arrangement

During the test, a horizontal force was applied to the top of the tower by a horizontal servo actuator. The vertical force was applied by the vertical jack. The moment applied a vertical force to the vertical jack above the tower, and an eccentric weight was applied to the top of the tower, as a moment. The specific installation situation is shown in Figure 4.

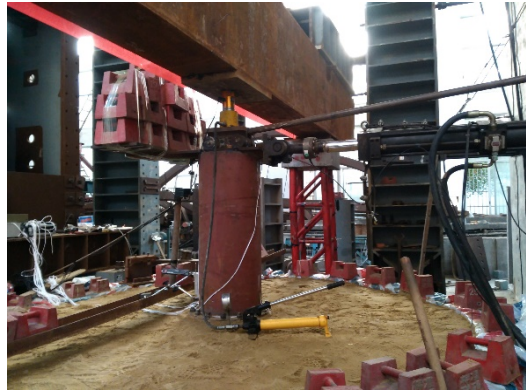

Figure 4. Loading device arrangement

To test the horizontal displacement and settlement displacement of the foundation under wind load, a vertical linear variable differential transformer and a horizon vertical linear variable differential with an effective stroke of 200 mm were arranged on the top of the foundation.

The strain is the main test object. In this study, multiple strain gauges were attached to the four parts of the test object. Firstly, to test the tower strain, a test section was arranged at a distance of 5 cm from the bottom of the tower, and the strain gauges were fixed at the outer barrel wall and an angle of 45°. Secondly, a vertical strain gauge was placed at a 45° angle in the middle of the base ring to test basic loop strain. Thirdly, to test the strain of the base plate, three lines were arranged on the base plate, and the strain gauges on the line were arranged along the wind direction, and the strain gauges passed through the base center. Fourthly, due to the rough surface of the nut connecting the base ring and the foundation, and the surface of the bolt was in contact with other parts, the strain gauge could not directly adhered. In this study, the bolt strain was indirectly measured by bonding steel sheets at the joint of the base ring and the tower and then by bonding the strain gauges on the steel sheet. However, considering that the steel sheet was located outside the bolt, and when the bolt is pressed, the steel sheet deformation and the bolt deformation are not linear. In this study, the steel sheet was only adhered on the windward side to test the tensile strain of the bolt.

The base pressure was tested by the earth pressure box. Due to the geometric symmetry of the model, the load is also symmetrically distributed. The direction from the wind load is the midline. Under the 1/2 basis, the distance from the base edge is 5 cm and the interval is 45°. One range is 300 MPa. Pressure box, and then placed a soil pressure box in the center of the base, a total of 6 earth pressure boxes, the arrangement of spinous processes is shown in Figure 5.

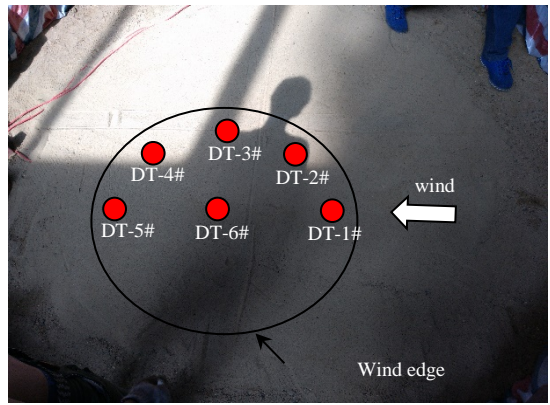

Figure 5. Earth pressure box arrangement

During the test, all test components were cleared after the application of vertical force and eccentric force was stabilized, that was, before the horizontal dynamic load was applied. The instrument only collects the displacement, settlement, strain and pressure caused by the horizontal dynamic load.

### 3.4 Analysis of test results

（1）Strain of the connecting bolt

As described above, because the surface shape and location of the bolt did not meet the conditions of use of the strain gauge, the test was to approximate the strain of the steel sheet to consider as bolt strain by attaching the steel sheet to the windward side of the joint. The test result is show in figure 6. The figure shows that the maximum strain is 28.3με, its corresponding axial stress is 5.94MPa. From the numerical point of view, the axial stress caused by the horizontal random wind load was small, but considering the random cyclic stress may cause fatigue damage of the steel, it cannot be completely ignored.

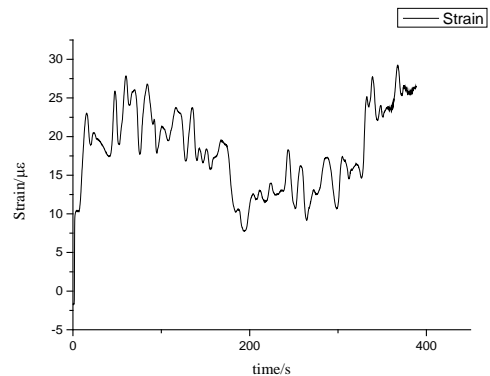

Figure 6. Windward side bolt vertical axial curve

（2）Vertical strain in the middle of the base ring

The vertical strain on the windward side and leeward side of the base ring is show in figure7 and figure8. The maximum tensile strain caused by the horizontal dynamic load on the basic ring is on the windward side, and the value is 17.2με. The maximum compressive strain caused by the horizontal dynamic load on the base ring is located on the wind side, and the maximum absolute value is 10.8με. According to the formula $\sigma = E\varepsilon$, the maximum tensile stress is 3.66MPa, the maximum  compressive stress is 2.298MPa.The value of stress is way less than steel yield strength 349.2MPa, so the possibility of damage to the base ring is small.

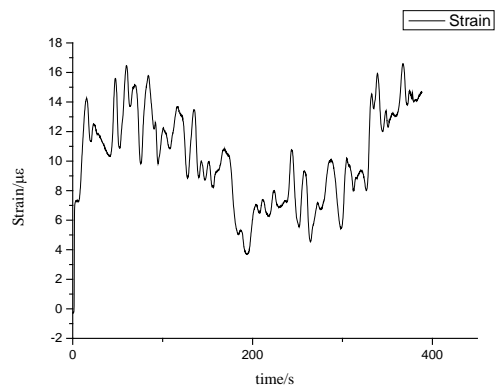

Figure 7. Vertical strain curve on the windward side of the base ring

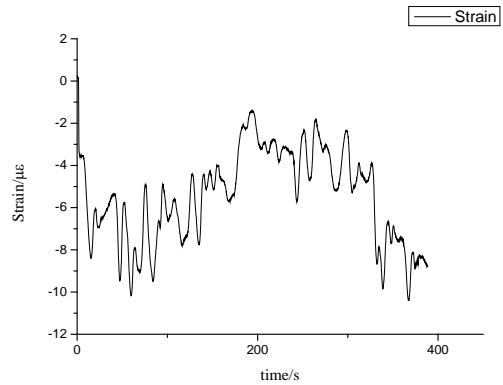

Figure 8. the vertical strain on the leeward side of the foundation ring

（3）Foundation displacement and settlement

The displacement-time curve and the settling-time curve are shown in Figure 9 and Figure 10. The horizontal displacement does not exceed 0.50mm, the settlement value does not exceed 0.14mm. Whether it is the displacement value or settlement value after multiplication by the proportionality factor of 10, it is less than the allowable range of 100mm.

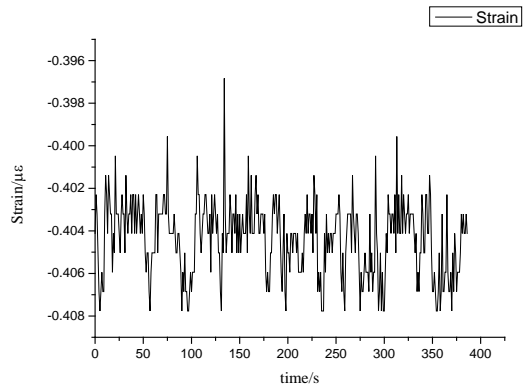

Figure 9. Foundation displacement-time curve

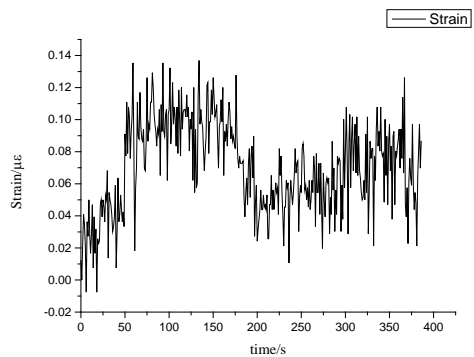

Figure 10. Foundation settlement-time curve

（4）Strain of the base plate

By analyzing the data collected by the strain gauges, it is found that the strain changes through the base center are the most severe. There are four strain gauges on this line, numbered JY-1#, JY-2#, JY-3# and JY-4# respectively. The maximum compressive strain is on the edge of the windward side foundation, and the value turns positive and gradually increases near the centerline of the substrate, and decreases at the position near the base edge of the leeward side. The specific change of the maximum substrate pressure on this line is as shown in Figure11. The value of maximum compressive is -10.6µε, it's way less than the ultimate compressive strain of concrete 3300µε. The value of tensile strain is 6.4µε, it's way less than the ultimate tensile strain of concrete100µε.

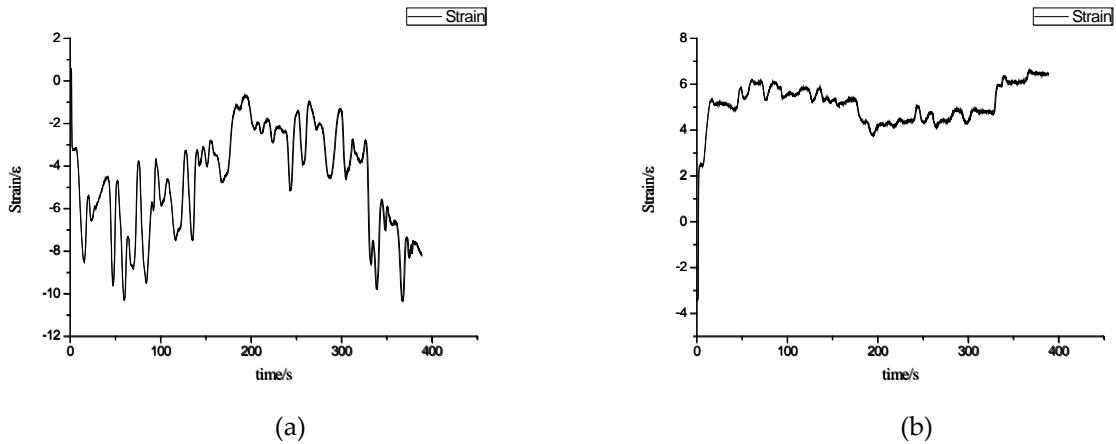

(a)                                          (b)

Figure 11. Base plate maximum strain-position curve. (a) JY-1#;(b) JY-3#.

（5）Foundation pressure

The change of the base pressure on the windward side and the wind side is shown in Figure 12. It can be seen from the figure that the change of the base pressure caused by the horizontal load is uneven, the base pressure caused on the windward side is small and some data is negative, the horizontal load makes the base pressure becomes smaller. The substrate pressure caused on the leeward side is positive and bigger than the base pressure caused on the windward side. The substrate load in the actual operation process is also unevenly distributed, but the base pressure is assumed to be evenly distributed in the specification, which has a certain influence on the design calculation.

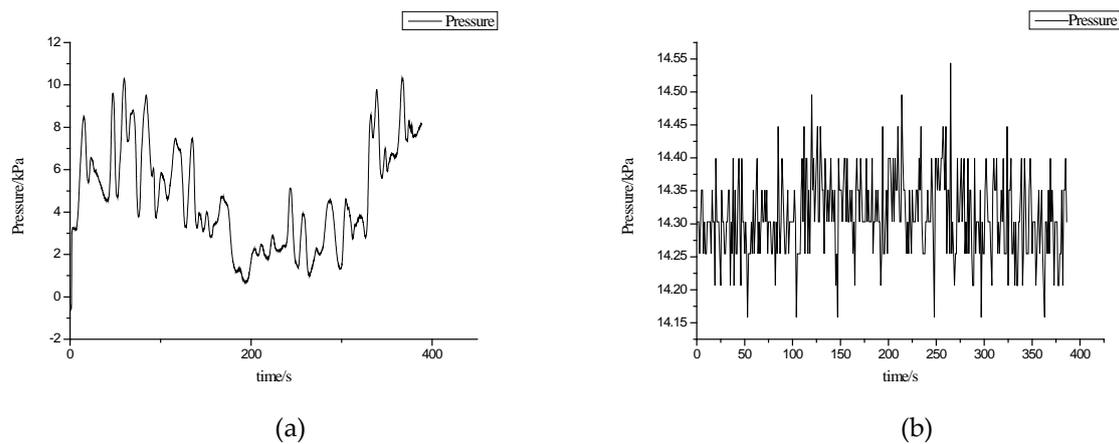

(a)                                          (b)

Figure 12. Foundation bottom pressure curve. (a) Windward side substrate pressure curve; (b) Gale side base pressure curve.

# 4. Numerical simulation

To study the influence of the uneven distribution of substrate pressure on wind power foundation damage, a numerical model was established to study the variation of foundation soil stress and foundation uneven settlement under horizontal wind load.

## 4.1 Establishment of the numerical model

According to the analysis in the previous section, the foundation calculation model of the wind turbine meets both the load symmetry condition and the geometric symmetry condition. Therefore, to simplify the calculation, 1/2 of the geometric model could be numerically modeled. The numerical model established is shown in Figure 13.

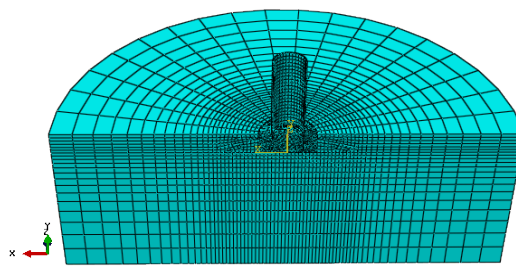

The infinite element method was used to solve the dynamic boundary problem, that is, the horizontally outward infinite unit CIN3D8 was used on the side of the foundation soil. The established model is show in Figure14. And the vertical force, horizontal resultant force and horizontal resultant force were applied to the top of the tower according to the simplified wind load and moment.

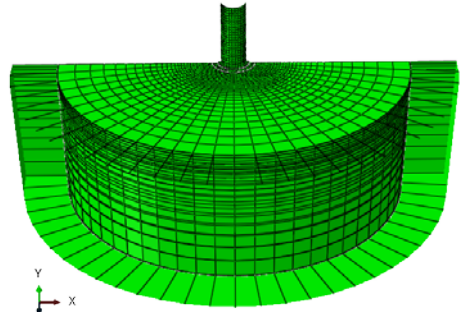

Figure 14. Infinite unit artificial boundary

## 4.2 Influence of wind load on the stress of foundation soil

By comparing the numerical simulation results with the model test basic ring and the basic strain, the foundation settlement and displacement, and the base pressure data, it was found that the data fluctuates in the same range, and the curve shape is the same. It was verified that this numerical model is suitable for this model structure. And the substrate pressure also has an uneven distribution that was different from the specification assumption. Taken this as a starting point, ABAQUS was used to simulate the distribution of the vertical stress of the foundation soil under the foundation edge when the wind turbine is subjected to the dynamic wind load and the static wind load respectively, and test the damage of the wind turbine foundation under strong wind load. Further explore the damage characteristics of the circular foundation of wind turbine with slope soil as the bearing layer, and study the influence of soil stress evolution caused by an uneven distribution of foundation pressure on foundation damage.

The test shows that compare with the static load, the stress of the foundation soil under the random dynamic wind load fluctuates significantly. The influence of the random wind load on the stress of the foundation soil could be attributed to the amplification and reduction of the stress of the foundation soil by the random wind load.

According to the finite element calculation results under dynamic and static wind loads, the vertical stress of the foundation soil changes with depth and width under the influence of random wind loads. The influence of the upper dynamic wind load on the foundation soil stress reaches 1.5 times the base radius, while the influence on the foundation soil stress is twice the base radius, but the core soil below the foundation is almost unaffected. Numerically, the maximum value of the dynamic load on the soil stress is below the edge of the windward side. The maximum amplification effect is higher than 1.64 times relative to the static load, and the soil stress below the leeward side is also affected by the dynamic load. The amplification effect is more than 1.20 times smaller than the windward side. It can be seen that the influence of dynamic wind load on the stress of foundation soil is not only significant in scope but also in numerical value.

## 4.3 Failure mode under strong wind loads

Based on the established numerical model, the horizontal wind load and its corresponding horizontal moment are gradually increased by the load of 0.1kN to test the damage of the foundation of the wind turbine. When the load increased to a certain value, the finite element calculation did not converge, and the last applied horizontal load is taken as the ultimate wind load of the foundation of the wind turbine. The equivalent plastic strain cloud diagram of the output wind turbine foundation failure is shown in Figure 15. It can be seen from the figure that the foundation of the wind turbine foundation is destroyed because the soil at the edge of the windward side enters the plastic zone, then local yielding occurs, and the backfill is filled on the foundation. In the uplift, the foundation on the windward side was lifted up significantly, and the foundation of the wind turbine was overturned.

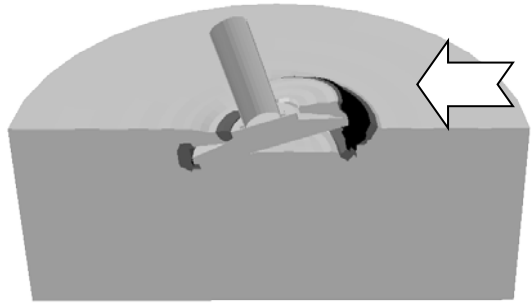

Figure 15. Wind turbine extension based on failure mode under strong wind loads

# 5. Conclusion

In this context, the circular expansion foundation of 2MW installed capacity was taken as the research object, and the 1:10 scale model test was established. Focus on the position of the foundation and the connection between the foundation and the tower, Bolt strain, base ring strain, base floor strain, base displacement and base pressure were tested in all directions. The failure mode of the enlarged concrete foundation with the slope soil as the bearing layer under strong wind load is studied. The main conclusions are as follows:

(1) Under the rated working conditions, bolt strain, base ring strain, base floor strain, base displacement and base pressure fluctuates within a reasonable range, but the influence of fatigue failure of the steel is still need considered.

(2) The wind-based dynamic response caused by random wind loads is the main reason cased the wind turbine damage, especially its influence on foundation soil stress is significant in scope and numerical value. The influence of the upper dynamic wind load on the foundation soil stress reaches 1.5 times the base radius, but the core soil below the foundation is almost unaffected. It is recommended that the soil in this range be reinforced in the construction project to ensure operational safety.

**Acknowledgments**

The authors are grateful for the valuable comments from the editors and the reviewer, which has substantially improved the quality of our work.

**Funding Statement**

This work was supported by the Open Fund of Key Laboratory of Earth and Rock Dam Failure Mechanism and Prevention and Control Technology of the Ministry of Water Resources (No.YK319008) and Doctoral Research Foundation of the University of South China (No.2012XQD01) .

**Data Accessibility**
**Competing Interests**
*We have no competing interests.*
**Authors' Contributions**
Peng Cheng conceived the article structure. Zijian Fan worte the paper